\documentclass[12pt]{amsart}
\usepackage{latexsym, amsmath, amssymb, amsthm}
\usepackage{epsfig}
\usepackage{amscd}
\usepackage{latexsym}
\usepackage{pstricks,pst-node,cite}

\newtheorem{thm}{Theorem}[section]

\newtheorem{cor}[thm]{Corollary}
\newtheorem{exa}[thm]{Example}

\newtheorem{defi}[thm]{Definition}
\theoremstyle{remark}

\newcommand{\F}{\mathcal{F}}

\newcommand{\BS}{\mathbb{S}}

\newcounter{fignum}
\ifx\blackandwhite\undefined
  \newrgbcolor{darkred}{.9 0 0}\newrgbcolor{emgreen}{0 .9 0}
  \newrgbcolor{pup}{.7  0 .9}\newrgbcolor{xxx}{.8 .8 .8}
  \newrgbcolor{lightgray}{.99 .99 0}\newrgbcolor{darkgray}{.8 .8 .3} \else
\fi
\psset{linecolor=blue}

\psset{linewidth=.7pt,dash=3pt 3pt,doublesep=.05,dotsize=1pt 5}
\SpecialCoor

\begin{document}

\author{Dongseok Kim}
\address{Department of Mathematics \\Kyonggi University
\\ Suwon, 443-760 Korea}
\email{dongseok@kgu.ac.kr}


\title[Links of the flat plumbing basket numbers $4$ or less]
{A classification of links of the flat plumbing basket numbers $4$ or less}

\begin{abstract}
Flat plumbing basket surfaces of links were introduced to study the geometry of the complement of the links.
In present article, we study links of the flat plumbing basket numbers $4$ or less using
a special presentation of the flat plumbing basket surfaces. We find a
complete classification theorem of links of the flat plumbing basket numbers $4$ or less.
\end{abstract}

\maketitle
\section{Introduction}

A \emph{link} $L$ is an embedding of $n$ copies of $\mathbb{S}^1$ in
$\mathbb{S}^3$. If the number of components of the link $L$ is $1$,
a link is called a \emph{knot}.
Throughout the article, we will assume all links are \emph{tame}
which means all links can be in a form of a finite union of line segments.
Two links are \emph{equivalent}
if there is an isotopy between them.
In the case of prime knots, this equivalence is the same as
the existence of an orientation preserving homeomorphism on $\mathbb{S}^3$,
which sends a knot to the other knot. Although the equivalent
class of a link $L$ is called a \emph{link type},
throughout the article, a link really means the equivalent class of the link.

A compact orientable surface $\F$ is called a \emph{Seifert surface} of a link $L$ if the boundary of $\F$ is isotopic to
$L$. The existence of such a surface was first proven by Seifert
using an algorithm on a diagram of $L$, this algorithm was named after him as
\emph{Seifert's algorithm}~\cite{Seifert:def}. A Seifert surface $\F_L$ of
an oriented link $L$ which is produced by applying Seifert's
algorithm to a link diagram and
is called a \emph{canonical Seifert surface}.

Some Seifert surfaces feature extra structures. Seifert
surfaces obtained by flat annuli plumbings are the main subjects of this
article. Even though higher dimensional plumbings can be defined,
but we will only concentrate on \emph{annuli plumbings}.
This is often called a \emph{Murasugi sum} and it
has been studied extensively for the fibreness of links and surfaces
\cite{Gabai:genera, Stallings:const}.

The definition of flat plumbing basket surfaces~\cite{Rudolph:plumbing} is very technical
and so it is difficult to handle but the work in~\cite{FHK:openbook} provided
a tangible equivalent definition of a flat plumbing basket surface using
an open book decomposition. In a recent article by Hirose and Nakashima~\cite{HN},
the flat plumbing basket numbers of knots of $9$ crossings or less
were studied using this definition. Using the results in~\cite{KKL:string} and in~\cite{Kim:flat},
Choi, Do and the author~\cite{CDK} are working on a new knot tabulation with respect
to the flat plumbing basket number of knots using Dowker-Thistlethwaite notation
and computer program {\tt{knotscape}}.
However, none of these methods can be directly applied for links with more than one components.

In present article, we study links of the flat plumbing basket numbers $4$ or less using
a special presentation of the flat plumbing basket surfaces. We find a
complete classification theorem of links of the flat plumbing basket numbers $4$ or less.

The outline of this paper is as follows. We first provide some preliminary definitions and results
in Section~\ref{prelim}.
We provide a special presentation of the flat plumbing basket surfaces with $n$-annuli where $n\le 4$.
Using this presentation, we prove complete classification theorem
of links of the flat plumbing basket numbers $4$ or less in Section~\ref{fpbs}.

\section{Preliminaries} \label{prelim}

The followings are the exact definitions of the flat
plumbing basket surface given by Rudolph~\cite{Rudolph:plumbing}.

Spaces, maps, etc., are piecewise smooth unless stated differently.
Let $M$ be an oriented manifold. ${-}M$ denotes $M$ with its orientation reversed
and when notation requires it, $+M$ denotes $M$. For a suitable subset
$S \subset M$, $N_M (S)$ denotes a closed regular neighborhood of $S$ in $(M, \partial M)$
where an ordered pair $(S, T)$ stands a condition $T \subset S$ and a
map between ordered pairs $f: (S, T) \rightarrow (U,V)$ is a
map $f: S \rightarrow U$ which requires to preserve
subsets so that $f(T) \subset V$.
For a suitable codimension-$1$ submanifold $S \subset M$ (resp.,
submanifold pair
$(S, \partial S) \subset (M, \partial M))$, a emph{collaring}
is an orientation-preserving embedding
$S \times [0, 1] \rightarrow M$ (resp., $(S, \partial S)
\times[0, 1] \rightarrow (M, \partial M))$ extending $id_S = id_{S \times\{0\}}$;
a \emph{collar} of $S$ in $M$ (resp., of $(S, \partial S)$
in $(M, \partial M)$) is the image col$_M (S)$
(resp., col$_{(M, \partial M)} (S, \partial S)$) of a collaring.
The push-off of $S$ determined by
a collaring of $S$ or $(S, \partial S)$, denoted by $S^+$,
is the image by the collaring of
$S \times \{1\}$ with the orientation of $S$;
let $S^{-} :=$ ${-}$ $S^{+}$ such that $S$ and $S^{-}$ are
oriented submanifolds of the boundary of col$_M (S)$.

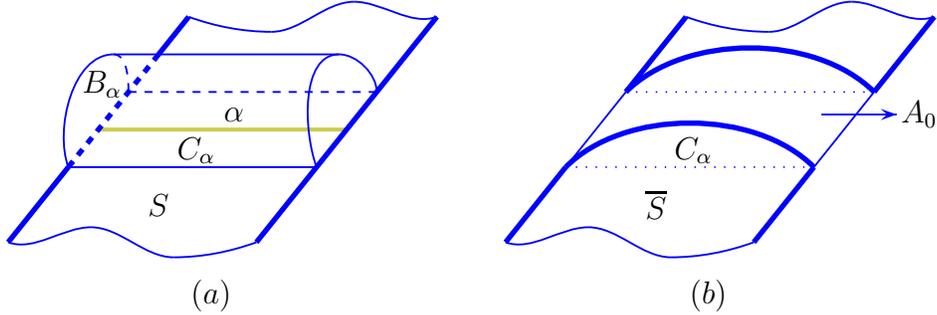
\begin{figure}
$$
\begin{pspicture}[shift=-1.6](-.2,-2.2)(6,1.7)
\psline[linewidth=2pt](0,-1.5)(.8,-.5)
\psline[linecolor=darkgray, linewidth=1.5pt](1.2,0)(4.5,0)
\psline[linestyle=dashed, linewidth=2pt](.8,-.5)(2,1)
\psline[linewidth=2pt](2.4,1.5)(2,1)
\pccurve[angleA=110,angleB=180](.8,-.5)(1.4,1)
\pccurve[linestyle=dashed, angleA=0,angleB=90](1.4,1)(1.6,.5)
\psline(1.4,1)(4.4,1)
\pccurve[angleA=110,angleB=180](4.1,-.5)(4.4,1)
\pccurve[angleA=0,angleB=90](4.4,1)(4.9,.5)
\psline[linewidth=2pt](3.3,-1.5)(5.7,1.5)
\psline(.8,-.5)(4.1,-.5) \psline[linestyle=dashed](1.6,.5)(4.9,.5)
\pccurve[angleA=200,angleB=0](5.7,1.5)(4.6,1.7)
\pccurve[angleA=180,angleB=0](4.6,1.7)(3.5,1.3)
\pccurve[angleA=180,angleB=-20](3.5,1.3)(2.4,1.5)
\pccurve[angleA=200,angleB=0](3.3,-1.5)(2.2,-1.7)
\pccurve[angleA=180,angleB=0](2.2,-1.7)(1.1,-1.3)
\pccurve[angleA=180,angleB=-20](1.1,-1.3)(0,-1.5)
\rput(3,.2){$\alpha$} \rput(2.5,-.28){$C_{\alpha}$}
\rput(2,-1){$S$} \rput(1.25,.6){$B_{\alpha}$}
\rput[t](2.7,-2){$(a)$}
\end{pspicture}\quad
\begin{pspicture}[shift=-1.6](-.2,-2.2)(6,1.7)
\psline[linewidth=2pt](0,-1.5)(.8,-.5)
\psline(.8,-.5)(1.6,.5) \psline[linewidth=2pt](1.6,.5)(2.4,1.5)
\pccurve[linewidth=2pt, angleA=45,angleB=135](1.6,.5)(4.9,.5)
\pccurve[linewidth=2pt, angleA=45,angleB=135](.8,-.5)(4.1,-.5)
\psline[linewidth=2pt](3.3,-1.5)(4.1,-.5)
\psline(4.1,-.5)(4.9,.5) \psline[linewidth=2pt](4.9,.5)(5.7,1.5)
\psline[linestyle=dotted](.8,-.5)(4.1,-.5)
\psline[linestyle=dotted](1.6,.5)(4.9,.5)
\pccurve[angleA=200,angleB=0](5.7,1.5)(4.6,1.7)
\pccurve[angleA=180,angleB=0](4.6,1.7)(3.5,1.3)
\pccurve[angleA=180,angleB=-20](3.5,1.3)(2.4,1.5)
\pccurve[angleA=200,angleB=0](3.3,-1.5)(2.2,-1.7)
\pccurve[angleA=180,angleB=0](2.2,-1.7)(1.1,-1.3)
\pccurve[angleA=180,angleB=-20](1.1,-1.3)(0,-1.5)
\psline(4.2,.2)(5.2,.2) \psline[arrowscale=1.5]{->}(5.18,.2)(5.2,.2)
\rput(2.5,-.28){$C_{\alpha}$} \rput(5.5,.2){$A_{0}$}
\rput(2,-1){$\overline{S}$}
\rput[t](2.7,-2){$(b)$}
\end{pspicture}
$$
\caption{$(a)$ A geometric shape of $\alpha, B_{\alpha}$
and $C_{\alpha}$ on a Seifert surface $S$
and $(b)$ a new Seifert surface $\overline{S}$ obtained
from $S$ by a top $A_0$ plumbing along the path
$\alpha$.}
\label{topfig}
\end{figure}

\begin{defi}(\cite{Rudolph:plumbing}) \label{topdef}
Let $\alpha$ be a proper arc on a Seifert surface $S$.
Let $C_{\alpha}$ be \emph{col}$_{(S, \partial S)} (\alpha, \partial \alpha)$ which is called
\emph{the gluing region}.
Let $B_{\alpha}$ be \emph{col}$_{(S, \partial S)} (\alpha, \partial \alpha)$ (so $B_{\alpha}$
is a $3$-cell in \emph{top}$(S)$, that is the positive normal to $S$
along $C_{\alpha}=S \cap B_{\alpha} \subset \partial B_{\alpha}$
points into $B_{\alpha}$) as depicted in Figure~\ref{topfig}.
Let $A_n\subset B_{\alpha}$ be an $n$-full twisted annulus such that
$A_n \cap \partial B_{\alpha}= C_{\alpha}$. Then \emph{top plumbing} on $S$ along a path $\alpha$
is the new surface $\overline{S}= S \cup A_n$ where $A_n, C_{\alpha}, B_{\alpha}$ satisfy
the previous conditions.
\end{defi}

\begin{defi}(\cite{Rudolph:plumbing}) \label{fpbdef}
A Seifert surface $S$ is a \emph{flat plumbing basket surface} if either it is
$2$-disc $D^2$ or it can be constructed by plumbing $A_0$ to a
flat plumbing basket surface $S_0$ along a proper arc $\alpha \subset D^2\subset S_0$.
We say that a link $L$ admits a \emph{flat plumbing basket
presentation} if there exists a flat plumbing basket $S$ such that
$\partial S$ is equivalent to $L$.
The \emph{flat plumbing basket number} of $L$, denoted
by $fpbk(L)$, is the minimal number of flat annuli to obtain a
flat plumbing basket surface of $L$.
\end{defi}

An alternative definition of the flat plumbing basket
surfaces is given in~\cite{FHK:openbook}
and it is very easy to follow. The \emph{trivial
open book decomposition} of $\mathbb{R}^3$
is a decomposition of $\mathbb{R}^3$ into the half planes
in the following form. In a cylindrical coordinate, it can be presented
$$ \mathbb{R}^3 = \bigcup_{\theta \in [0, 2\pi)}
\{(r, \theta, z) | r \ge 0, z \in \mathbb{R} \}$$
where $\{(r, \theta, z) | r \ge 0, z \in \mathbb{R} \}$
is called a \emph{page} for $\theta \in [0, 2\pi)$.
Let $\mathcal{O}$ be the \emph{trivial open book decomposition} of the
$3$-sphere $\BS^3$ which is obtained from the trivial
open book decomposition of $\mathbb{R}^3$
by the one point compactification. A Seifert surface
is said to be a flat plumbing basket surface
if it consists of a single page of $\mathcal{O}$ as a $2$-disc $D^2$ and
finitely many bands which are
embedded in distinct pages~\cite{FHK:openbook}.

\begin{figure}
$$
\begin{pspicture}[shift=-1.2](0,-1.8)(9.5,3.1)
\psline[linecolor=pup, linestyle=dashed](.05,.05)(.05,1.95)
\pccurve[linecolor=pup, linestyle=dashed, angleA=90,angleB=180](.05,1.95)(1.05,2.95)
\psline[linecolor=pup, linestyle=dashed](1.05,2.95)(8.15,2.95)
\pccurve[linecolor=pup, linestyle=dashed, angleA=0,angleB=90](8.15,2.95)(9.15,2.05)
\psline[linecolor=pup, linestyle=dashed](9.15,2.05)(9.15,.05)
\pccurve[linecolor=pup, linestyle=dashed, angleA=-90,angleB=0](9.15,.05)(8.15,-.95)
\psline[linecolor=pup, linestyle=dashed](8.15,-.95)(2.95,-.95)
\pccurve[linecolor=pup, linestyle=dashed, angleA=180,angleB=-45](2.95,-.95)(2.65,-.55)
\pccurve[linecolor=pup, linestyle=dashed, angleA=135,angleB=0](2.45,-.35)(2.15,.05)
\pccurve[linecolor=pup, linestyle=dashed, angleA=180,angleB=0](2.15,.05)(1.35,1.05)
\pccurve[linecolor=pup, linestyle=dashed, angleA=180,angleB=180](1.35,1.05)(1.35,1.95)
\psline[linecolor=pup, linestyle=dashed](1.35,1.95)(7.85,1.95)
\pccurve[linecolor=pup, linestyle=dashed, angleA=0,angleB=0](7.85,1.95)(7.85,1.05)
\psline[linecolor=pup, linestyle=dashed](7.85,1.05)(2.05,1.05)
\pccurve[linecolor=pup, linestyle=dashed, angleA=180,angleB=45](2.05,1.000005)(1.75,.65)
\pccurve[linecolor=pup, linestyle=dashed, angleA=-135,angleB=0](1.55,.45)(1.25,.05)
\psline[linecolor=pup, linestyle=dashed](1.25,.05)(1.05,.05)
\pccurve[linecolor=pup, linestyle=dashed, angleA=180,angleB=-90](1.05,.05)(.55,.55)
\psline[linecolor=pup, linestyle=dashed](.55,.55)(.55,2.05)
\pccurve[linecolor=pup, linestyle=dashed, angleA=90,angleB=180](.55,2.05)(1.05,2.45)
\psline[linecolor=pup, linestyle=dashed](1.05,2.45)(8.15,2.45)
\pccurve[linecolor=pup, linestyle=dashed, angleA=0,angleB=90](8.15,2.45)(8.65,1.95)
\psline[linecolor=pup, linestyle=dashed](8.65,1.95)(8.65,.45)
\pccurve[linecolor=pup, linestyle=dashed, angleA=-90,angleB=0](8.65,.45)(8.15,.05)
\psline[linecolor=pup, linestyle=dashed](8.15,.05)(2.85,.05)
\pccurve[linecolor=pup, linestyle=dashed, angleA=180,angleB=0](2.85,.05)(2.05,-.95)
\psline[linecolor=pup, linestyle=dashed](2.05,-.95)(1.05,-.95)
\pccurve[linecolor=pup, linestyle=dashed, angleA=180,angleB=-90](1.05,-.95)(.05,.05)
\psline(0,0)(0,2)
\pccurve[angleA=90,angleB=180](0,2)(1,3)
\psline(1,3)(8.2,3)
\pccurve[angleA=0,angleB=90](8.2,3)(9.2,2)
\psline(9.2,2)(9.2,0)
\pccurve[angleA=-90,angleB=0](9.2,0)(8.2,-1)
\psline(8.2,-1)(6.9,-1)
\pccurve[angleA=180,angleB=-45](6.9,-1)(6.6,-.6)
\pccurve[angleA=135,angleB=0](6.4,-.4)(6.1,0)
\psline(6.1,0)(5.3,0)
\pccurve[angleA=180,angleB=-45](5.3,0)(5,.4)
\pccurve[angleA=135,angleB=0](4.8,.6)(4.5,1)
\pccurve[angleA=180,angleB=0](4.5,1)(3.7,0)
\pccurve[angleA=180,angleB=-45](3.7,0)(3.4,.4)
\pccurve[angleA=135,angleB=0](3.2,.6)(2.9,1)
\psline(2.9,1)(2.1,1)
\pccurve[angleA=180,angleB=45](2.1,1)(1.8,.6)
\pccurve[angleA=-135,angleB=0](1.6,.4)(1.3,0)
\psline(1.3,0)(1,0)
\pccurve[angleA=180,angleB=-90](1,0)(.5,.5)
\psline(.5,.5)(.5,2)
\pccurve[angleA=90,angleB=180](.5,2)(1,2.5)
\psline(1,2.5)(8.2,2.5)
\pccurve[angleA=0,angleB=90](8.2,2.5)(8.7,2)
\psline(8.7,2)(8.7,.5)
\pccurve[angleA=-90,angleB=0](8.7,.5)(8.2,0)
\psline(8.2,0)(6.9,0)
\pccurve[angleA=180,angleB=0](6.9,0)(6.1,-1)
\psline(6.1,-1)(2.9,-1)
\pccurve[angleA=180,angleB=-45](2.9,-1)(2.6,-.6)
\pccurve[angleA=135,angleB=0](2.4,-.4)(2.1,0)
\pccurve[angleA=180,angleB=0](2.1,0)(1.3,1)
\pccurve[angleA=180,angleB=180](1.3,1)(1.3,2)
\psline(1.3,2)(7.9,2)
\pccurve[angleA=0,angleB=0](7.9,2)(7.9,1)
\psline(7.9,1)(5.3,1)
\pccurve[angleA=180,angleB=0](5.3,1)(4.5,0)
\pccurve[angleA=180,angleB=-45](4.5,0)(4.2,.4)
\pccurve[angleA=135,angleB=0](4,.6)(3.7,1)
\pccurve[angleA=180,angleB=0](3.7,1)(2.9,0)
\pccurve[angleA=180,angleB=0](2.9,0)(2.1,-1)
\psline(2.1,-1)(1,-1)
\pccurve[angleA=180,angleB=-90](1,-1)(0,0)
\rput(4.6,-1.6){{$(a)$}}
\end{pspicture}$$
$$
\begin{pspicture}[shift=-1.2](0,-2)(9.5,3.1)
\psline[linecolor=darkred, linestyle=dashed](2.9,1)(6.9,1)
\psline[linecolor=darkred, linestyle=dashed](2.9,0)(6.9,0)
\psline[linecolor=darkred, linestyle=dashed](2.9,-1)(6.9,-1)
\psline(0,0)(0,2)
\pccurve[angleA=90,angleB=180](0,2)(1,3)
\psline(1,3)(8.2,3)
\pccurve[angleA=0,angleB=90](8.2,3)(9.2,2)
\psline(9.2,2)(9.2,0)
\pccurve[angleA=-90,angleB=0](9.2,0)(8.2,-1)
\psline(8.2,-1)(6.9,-1)
\pccurve[angleA=180,angleB=0](6.9,-1)(6.1,0)
\psline(6.1,0)(5.3,0)
\pccurve[angleA=180,angleB=-45](5.3,0)(5,.4)
\pccurve[angleA=135,angleB=0](4.8,.6)(4.5,1)
\pccurve[angleA=180,angleB=0](4.5,1)(3.7,0)
\pccurve[angleA=180,angleB=-45](3.7,0)(3.4,.4)
\pccurve[angleA=135,angleB=0](3.2,.6)(2.9,1)
\psline(2.9,1)(2.1,1)
\pccurve[angleA=180,angleB=45](2.1,1)(1.8,.6)
\pccurve[angleA=-135,angleB=0](1.6,.4)(1.3,0)
\psline(1.3,0)(1,0)
\pccurve[angleA=180,angleB=-90](1,0)(.5,.5)
\psline(.5,.5)(.5,2)
\pccurve[angleA=90,angleB=180](.5,2)(1,2.5)
\psline(1,2.5)(8.2,2.5)
\pccurve[angleA=0,angleB=90](8.2,2.5)(8.7,2)
\psline(8.7,2)(8.7,.5)
\pccurve[angleA=-90,angleB=0](8.7,.5)(8.2,0)
\psline(8.2,0)(6.9,0)
\pccurve[angleA=180,angleB=45](6.9,0)(6.6,-.4)
\pccurve[angleA=-135,angleB=0](6.4,-.6)(6.1,-1)
\psline(6.1,-1)(2.9,-1)
\pccurve[angleA=180,angleB=-45](2.9,-1)(2.6,-.6)
\pccurve[angleA=135,angleB=0](2.4,-.4)(2.1,0)
\pccurve[angleA=180,angleB=0](2.1,0)(1.3,1)
\pccurve[angleA=180,angleB=180](1.3,1)(1.3,2)
\psline(1.3,2)(7.9,2)
\pccurve[angleA=0,angleB=0](7.9,2)(7.9,1)
\psline(7.9,1)(5.3,1)
\pccurve[angleA=180,angleB=0](5.3,1)(4.5,0)
\pccurve[angleA=180,angleB=-45](4.5,0)(4.2,.4)
\pccurve[angleA=135,angleB=0](4,.6)(3.7,1)
\pccurve[angleA=180,angleB=0](3.7,1)(2.9,0)
\pccurve[angleA=180,angleB=0](2.9,0)(2.1,-1)
\psline(2.1,-1)(1,-1)
\pccurve[angleA=180,angleB=-90](1,-1)(0,0)
\pscircle[linecolor=blue, fillcolor=darkred, fillstyle=solid](1.3,-1){.1}
\psline[arrowscale=1.5]{->}(1.8,-1)(1.9,-1)
\psarc[doubleline=true](6.35,0){.75}{0}{180}
\psarc[doubleline=true](6.75,0){.75}{0}{180}
\rput(3.3,-.3){{$1$}}
\rput(4.1,-.3){{$2$}}
\rput(4.9,-.3){{$3$}}
\rput(5.5,-.3){{$4$}}
\rput(5.9,-.3){{$5$}}
\rput(6.5,.3){{$6$}}
\rput(7.1,-.3){{$7$}}
\rput(7.5,-.3){{$8$}}
\rput(3.3,1.3){{$9$}}
\rput(4.1,1.3){{$10$}}
\rput(4.9,1.3){{$11$}}
\rput(6.5,-1.3){{$12$}}
\rput(4.6,-1.8){{$(b)$}}
\end{pspicture}$$
$$
\begin{pspicture}[shift=-.8](-.7,-2.4)(6.2,2)
\psarc[doubleline=true](4,-.5){1.5}{-5}{185}
\psarc[doubleline=true](2.75,-.5){.75}{-5}{185}
\psarc[doubleline=true](2.25,-.5){.75}{-5}{185}
\psarc[doubleline=true](3,-.5){2}{-5}{185}
\psarc[doubleline=true](2.5,-.5){2}{-5}{185}
\psarc[doubleline=true](2,-.5){2}{-5}{185}
\psframe[linecolor=lightgray,fillstyle=solid,fillcolor=lightgray](-.5,-1.5)(6,-.5)
\psline(-.03,-.5)(-.5,-.5)(-.5,-1.5)(6,-1.5)(6,-.5)(4.53,-.5)
\psline(.03,-.5)(.47,-.5) \psline(.53,-.5)(.97,-.5) \psline(1.03,-.5)(1.47,-.5)
\psline(1.53,-.5)(1.97,-.5) \psline(2.03,-.5)(2.47,-.5) \psline(2.53,-.5)(2.97,-.5)
\psline(3.03,-.5)(3.47,-.5) \psline(3.53,-.5)(3.97,-.5) \psline(4.03,-.5)(4.47,-.5)
\psline(4.53,-.5)(4.97,-.5) \psline(5.03,-.5)(5.47,-.5)
\rput(0,-.75){{$1$}} \rput(.5,-.75){{$2$}}
\rput(1,-.75){{$3$}} \rput(1.5,-.75){{$4$}}
\rput(2,-.75){{$5$}} \rput(2.5,-.75){{$6$}}
\rput(3,-.75){{$7$}} \rput(3.5,-.75){{$8$}}
\rput(4,-.75){{$9$}} \rput(4.5,-.75){{$10$}}
\rput(5,-.75){{$11$}} \rput(5.5,-.75){{$12$}}
\rput(2.75,-1.2){{$\mathcal{D}$}}
\rput(2.75,-2.1){{$(c)$}}
\end{pspicture}
$$
\caption{$(a)$ The knot $5_2$ as a closed braid, $(b)$ Seifert surface of $5_2$ in order to apply the algorithm in~\cite{FHK:openbook},
$(c)$ a flat plumbing basket surface of $5_2$.} \label{52complete}
\end{figure}

In~\cite{FHK:openbook}, it is
shown that every link admits a flat plumbing basket representation by setting up the link
as a special closed braid form. So
we can define the \emph{flat plumbing basket number} of $L$, denoted
by $fpbk(L)$, to be the minimal number of flat annuli to obtain a
flat plumbing basket surface of $L$.
The author proved that every link $L$ admits a flat plumbing basket representation
by modifying the Seifert surface $S_L$ of the link $L$ to have a property that the Seifert graph
$\Gamma(S_L)$ holds the property described in~\cite[Theorem 3.3]{Kim:flat}.

Furthermore, an algorithm provided in~\cite{FHK:openbook} and in~\cite{CDK},
every links admit a flat plumbing basket presentation
$(a_1, a_2, \ldots, a_{2n})$ where $a_i \in \{1, 2, \ldots, n\}$
and each $i\in  \{1, 2, \ldots, n\}$ appears exactly twice.
In stead of explaining the exact algorithm, let us
provide a concrete example as follows.

\begin{exa}
The knot $5_2$ has a flat plumbing basket presentation
$(1,2,3,4,5,6,4,5,1$, $2$, $3$, $6)$ as
illustrated in Figure~\ref{52complete}.
\begin{proof}
From the braid representative $\sigma_2\sigma_1^{-1} (\sigma_2^{-1})^3 \sigma_1^{-1}$ of the knot $5_2$,
we choose the $2$-disc $\mathcal{D}$ the union of three disc bounded by three Seifert circle and two
half twisted bands represented by $\sigma_2\sigma_1^{-1}$ as depicted as the dashed purple line in
Figure~\ref{52complete} $(a)$. To make flat plumbing
we first change the crossing presented by $\sigma_1^{-1}$ by adding extra two flat annuli
as shown in Figure~\ref{52complete} $(b)$. Fix a starting point and an orientation coming from the braid
as indicated the red ball and arrow in Figure~\ref{52complete} $(b)$.
Put numbering for each bands when we move around the $2$-disc $\mathcal{D}$ from the starting point
in the direction indicated as given in Figure~\ref{52complete} $(b)$. Isotop the $2$-disc $\mathcal{D}$
into the standard $2$-disc $\mathcal{D}$ as illustrated in Figure~\ref{52complete} $(c)$.
Now we are ready to find a flat plumbing basket presentation of the flat plumbing basket surface
which is in the position of the trivial open book decomposition of $\BS^3$.
From the top annulus, a pair of two points in the boundary of the annulus will receive $1$ and so on.
In Figure~\ref{52complete} $(c)$, there is no distinction of order between two groups annuli
connecting $\{(1,7), (2,8),(3,9)\}$ and $\{(4,7),(5,8)\}$. If we consider $\{(1,7), (2,8),(3,9)\}$
are in front of $\{(4,7),(5,8)\}$, we get a flat plumbing basket presentation
$(1$, $2$, $3,4,5,6,4,5,1$, $2$, $3$, $6)$. Otherwise, one may get
$(3$, $4$, $5,1,2,6,1,2,3$, $4$, $5$, $6)$.
\end{proof}
\end{exa}

\begin{figure}
$$
\begin{pspicture}[shift=-1.2](-.7,-1.8)(4.2,1.2)
\psarc[doubleline=true](2.5,0){1}{-5}{185}
\psarc[doubleline=true](2,0){1}{-5}{185}
\psarc[doubleline=true](1.5,0){1}{-5}{185}
\psarc[doubleline=true](1,0){1}{-5}{185}
\psframe[linecolor=lightgray,fillstyle=solid,fillcolor=lightgray](-.5,-1)(4,0)
\psline(-.03,0)(-.5,0)(-.5,-1)(4,-1)(4,0)(3.53,0)
\psline(.03,0)(.47,0) \psline(.53,0)(.97,0) \psline(1.03,0)(1.47,0)
\psline(1.53,0)(1.97,0) \psline(2.03,0)(2.47,0) \psline(2.53,0)(2.97,0)
\psline(3.03,0)(3.47,0)
\rput(1.75,-.5){{$\mathcal{D}$}}
\rput(1.75,-1.7){{$(a)$}}
\end{pspicture} \quad
\begin{pspicture}[shift=-1.2](-.7,-1.8)(4.2,1.2)
\psarc[doubleline=true](2,0){1}{-5}{185}
\psarc[doubleline=true](2.5,0){1}{-5}{185}
\psarc[doubleline=true](1.5,0){1}{-5}{185}
\psarc[doubleline=true](1,0){1}{-5}{185}
\psframe[linecolor=lightgray,fillstyle=solid,fillcolor=lightgray](-.5,-1)(4,0)
\psline(-.03,0)(-.5,0)(-.5,-1)(4,-1)(4,0)(3.53,0)
\psline(.03,0)(.47,0) \psline(.53,0)(.97,0) \psline(1.03,0)(1.47,0)
\psline(1.53,0)(1.97,0) \psline(2.03,0)(2.47,0) \psline(2.53,0)(2.97,0)
\psline(3.03,0)(3.47,0)
\rput(1.75,-.5){{$\mathcal{D}$}}
\rput(1.75,-1.7){{$(b)$}}
\end{pspicture}
$$
\caption{Flat plumbing basket surfaces of $(a)$ the trefoil knot and $(b)$ the figure eight knot.} \label{figure34}
\end{figure}

\begin{exa} (\cite{CDK})
Flat plumbing basket numbers of the trefoil knot and the figure eight knot are $4$.
\begin{proof}
Flat plumbing basket surfaces of the trefoil knot and the figure eight knot with four annuli are depicted in Figure~\ref{figure34}.
By Theorem~\cite[Theorem 3.3]{CDK}, the flat plumbing basket number is bigger than or equal to the flat band index of the link.
The flat band index of the trefoil knot and the figure eight knot are $4$~\cite{KKL:string}. It complete the proof.
\end{proof}
\end{exa}

One may notice that these flat plumbing basket surfaces of three knots $3_1$, $4_1$ and $5_2$
have a common property that all of numbers in the set $\{ 1, 2, \ldots, n\}$ appear in the first half and
the last half of the flat plumbing basket presentation. Thus, we may
rewrite it as a \emph{permutations presentation} $(\sigma : \mu)$,
the first permutation $\sigma$ presents the order of
the annuli in the flat plumbing basket surface are connected, which is called the
\emph{connection permutation} and the second permutation $\mu$
presents the order of annuli from the top to the bottom
which is called the \emph{order permutation}. For example,
the flat plumbing basket surfaces of knot $5_2$ in Figure~\ref{52complete}
has permutations presentation $(3,4,5,1,2,6 : 1,2,3,4,5,6)$
and $3_1$ has $(1,2,3,4:1,2,3,4)$ while $4_1$ has
$(1,2,3,4: 1,2,4,3)$.

\section{Classification} \label{fpbs}

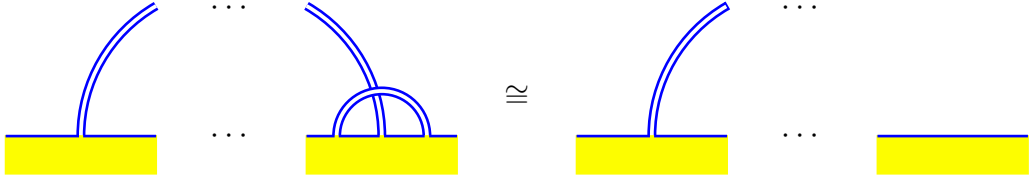
\begin{figure}
$$
\begin{pspicture}[shift=-1.2](-3.1,-.8)(3.1,2.2)
\psarc[doubleline=true, linewidth=1pt](0,0){2}{-10}{60}
\psarc[doubleline=true, linewidth=1pt](0,0){2}{120}{190}
\psarc[doubleline=true, linewidth=1pt](2,0){.6}{-10}{190}
\pspolygon[linecolor=lightgray, fillcolor=lightgray, fillstyle=solid](-1,0)(-3,0)(-3,-.5)(-1,-.5)(-1,0)
\pspolygon[linecolor=lightgray, fillcolor=lightgray, fillstyle=solid](1,0)(3,0)(3,-.5)(1,-.5)(1,0)
\psline[linewidth=1pt](-3,0)(-2.025,0)
\psline[linewidth=1pt](-1.97,0)(-1,0)
\psline[linewidth=1pt](1,0)(1.373,0)
\psline[linewidth=1pt](1.428,0)(1.972,0)
\psline[linewidth=1pt](2.03,0)(2.57,0)
\psline[linewidth=1pt](2.63,0)(3,0)
\rput(0,1.7){{$\cdots$}}
\rput(0,0){{$\cdots$}}
\end{pspicture}
\quad \cong \quad
\begin{pspicture}[shift=-1.2](-3.1,-.8)(3.1,2.2)
\psarc[doubleline=true, linewidth=1pt](0,0){2}{120}{190}
\psline(1.94;120)(2.06;120)
\pspolygon[linecolor=lightgray, fillcolor=lightgray, fillstyle=solid](1,0)(3,0)(3,-.5)(1,-.5)(1,0)
\pspolygon[linecolor=lightgray, fillcolor=lightgray, fillstyle=solid](-1,0)(-3,0)(-3,-.5)(-1,-.5)(-1,0)
\psline[linewidth=1pt](-3,0)(-2.025,0)
\psline[linewidth=1pt](-1.97,0)(-1,0)
\psline[linewidth=1pt](1,0)(3,0)
\rput(0,1.7){{$\cdots$}}
\rput(0,0){{$\cdots$}}
\end{pspicture}
$$
\caption{A fundamental move which decreases the flat plumbing basket number by $2$.} \label{twobridge}
\end{figure}

Since the flat plumbing basket number is defined to be the minimal number of flat annuli to obtain a
flat plumbing basket surface of $L$. If some annuli in flat plumbing basket presentations
can be removed, we should not consider such presentation from the beginning.
If there is a part of the form $iji$ in a flat plumbing basket presentation,
two annuli presented by $i$ and $j$ can be removed by a simple isotopy as illustrated in Figure~\ref{twobridge}
and such a flat plumbing basket presentation is said to be \emph{reducible}.
If a permutations presentation $(a_1, a_2, \ldots, a_n : b_1, b_2 , \ldots, b_n)$
admits such a removal of annuli, we say that the permutations presentation is \emph{reducible}.

First we show that any flat plumbing basket surface with $4$ or less
admits a permutations presentation as follows.

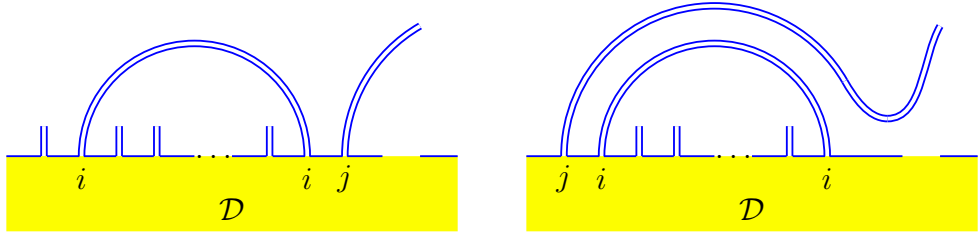
\begin{figure}
$$
\begin{pspicture}[shift=-1.2](-2.7,-1.2)(3.8,3.2)
\psarc[doubleline=true](0,0){1.5}{-5}{185}
\psarc[doubleline=true](4,0){2}{120}{185}
\psline[doubleline=true](-1,0)(-1,.4)
\psline[doubleline=true](-.5,0)(-.5,.4)
\psline[doubleline=true](1,0)(1,.4)
\psline[doubleline=true](-2,0)(-2,.4)
\psframe[linecolor=lightgray,fillstyle=solid,fillcolor=lightgray](-2.5,-1)(3.5,0)
\psline(-1.53,0)(-1.96,0)
\psline(-2.04,0)(-2.5,0)
\psline(-1.47,0)(-1.04,0) \psline(-.96,0)(-.54,0) \psline(-.46,0)(0,0)
\psline(.5,0)(.96,0) \psline(1.04,0)(1.47,0)
\psline(1.53,0)(1.97,0) \psline(2.03,0)(2.5,0) \psline(3,0)(3.5,0)
\rput(-1.5,-.3){{$i$}}
\rput(1.5,-.3){{$i$}}
\rput(2,-.3){{$j$}}
\rput(.25,0){{$\ldots$}}
\rput(.5,-.7){{$\mathcal{D}$}}
\end{pspicture}
\quad
\begin{pspicture}[shift=-1.2](-2.7,-1.2)(3.8,3.2)
\psarc[doubleline=true](0,0){2}{30}{185}
\psarc[doubleline=true](0,0){1.5}{-5}{185}
\psline[doubleline=true](-1,0)(-1,.4)
\psline[doubleline=true](-.5,0)(-.5,.4)
\psline[doubleline=true](1,0)(1,.4)
\pccurve[doubleline=true, angleA=-60,angleB=180](2;30)(2.3,.5)
\pccurve[doubleline=true, angleA=0,angleB=-120](2.3,.5)(3,1.732)
\psframe[linecolor=lightgray,fillstyle=solid,fillcolor=lightgray](-2.5,-1)(3.5,0)
\psline(-1.53,0)(-1.96,0)
\psline(-2.04,0)(-2.5,0)
\psline(-1.47,0)(-1.04,0) \psline(-.96,0)(-.54,0) \psline(-.46,0)(0,0)
\psline(.5,0)(.96,0) \psline(1.04,0)(1.47,0)
\psline(1.53,0)(2.5,0) \psline(3,0)(3.5,0)
\rput(-1.5,-.3){{$i$}}
\rput(1.5,-.3){{$i$}}
\rput(-2,-.3){{$j$}}
\rput(.25,0){{$\ldots$}}
\rput(.5,-.7){{$\mathcal{D}$}}
\end{pspicture}
$$
\caption{A handle slide the band which is labeled by $j$ along the
band labeled by $i$.} \label{slide}
\end{figure}

\begin{thm} \label{twopermutation}
The flat plumbing basket surfaces with $4$ or less annuli can be
presented by a permutations presentation.
\begin{proof}
For $n \le 1$, is obvious. For $n=2$, there are six flat plumbing basket presentations
$(i,i,j,j)$ $(i,j,i,j)$ where $i,j=1,2$.
The first one $(i,i,j,j)$ can be changed to $(i,j,j,i)$ by a
slide of annulus presented by $i$ along the annulus presented by $j$ as illustrated in Figure~\ref{slide}.
The second $(i,j,i,j)$ flat plumbing basket presentations
is already written as a permutations presentation without changing anything.
For $n=3$, if two numbers are in the first half, there are three possibilities :
either $i,i,j$, $j,i,i$ or $i,j,i$ where $i, j=1, 2, 3$.
But the last one is reducible to a flat plumbing basket presentation $(1,1)$
by the move described in Figure~\ref{twobridge}.
For the first one, it is either $(i,i,j,k,k,j)$ or $(i,i,j,j,k,k)$ for some $k\in \{1,2,3\}$.
By sliding the annulus presented by $i$ along the
annulus presented by $j$ for $(i,i,j,k,k,j)$, we get $(i,j,k,k,j,i)$ which admits a permutations presentation.
Cyclically rotating $(i,i,j,j,k,k)$ by changing the starting point, we get $(i,j,j,k,k,i)$ we reduce to a previous case.
For the second, it is either $(j,i,i,j,k,k)$ or $(j,i,i,k,k,j)$ and each of them can be
transformed to one admits a permutations presentation by a handle slide.

For $n=4$, now we divide cases depending on the number of the same letter in the first half, say $m$.
If $m=0$, we are done. If $m=1$ and the flat plumbing basket presentation is not reducible,
it is either $(i,i,j,k,*,*,*,*)$, $(i, j, j, k, *,*,*,*)$, $(i, j, k, k, *,*,*,*)$ or
$(i, j, k, i, *,*,*,*)$ where $i,j,k=1,2,3,$.
For $(i,i,j,k,*,*,*,*)$, the fifth components must be either $k$ or $l$ for some $l\in\{1,2,3,4\}$.
If it is $l$, by sliding the annulus presented by $i$ along the
annulus presented by $j$ for $(i,i,j,k,l,*,*,*)$, we have $(i,j,k,l,*,*,*,*)$ which admits a permutations presentation.
If it is $k$, there are two possibilities as follows and each can be changed to
one which admits a permutations presentation by sliding twice,

$$
\begin{matrix}
(1)& (i,i,j,k,k,l,l,j) & \rightarrow & (i,j,k,k,l,l,j,i) & \rightarrow & (i,j,k,l,l,k,j,i) \\
(2)& (i,i,j,k,k,j,l,l) & \rightarrow & (i,j,k,k,j,i,l,l) & \rightarrow & (l,i,j,k,k,j,i,l)
\end{matrix}
$$

For $(i, j, j, k, *,*,*,*)$, we only need to look at $(i,j,j,k,i,l,l,k)$ which can be changed
to $(i,j,k,i,l,l,k,j)$ and then to $(l,i,j,k,i,l,k,j)$ by two handle slide.
For $(i, j, k, k, *,*,*,*)$, there are six
possibilities but all can be changed to one which admits a permutations presentation by sliding
the annulus presented by $l$ along one annulus which is the right before $l$.
For $(i, j, k, i, *,*,*,*)$, we look the second half. By the previous argument, the only cases we have
to take care of is $(i, j, k, i, l, j, k, l)$ or $(i, j, k, i, l, k, j, l)$. But both cases
can be changed to one which admits a permutations presentation by sliding
the annulus presented by $i$ along the annulus presented by $l$.

If $m=2$, there are four possibilities, $(i,j,j,i,k,l,l,k)$, $(i,j,j,i,k,k,l,l)$,
$(i,i,j,j,k,l,l,k)$ and $(i,i,j,j,k,k,l,l)$. But any sliding of an annulus along an adjacent annulus
will decrease $m$ by $1$, it returns to the previous cases. Therefore, it completes the proof of theorem.
\end{proof}
\end{thm}

Now, we are set to prove the main theorem which completely classifies the links
of the flat plumbing basket numbers $4$ or less.

\begin{table}
\begin{tabular}{|c|c|c|c|c|}\hline
\rm{Name ~of ~link} & $\begin{matrix} \rm{Connection} \\ \rm{Permutation} \end{matrix}$
 & $\begin{matrix} \rm{Order} \\ \rm{Permutation} \end{matrix}$  \\
\hline $3_1$ & $1234$  & $\begin{matrix} 1234, 1432, 2143, 2341,\\ 3214, 3412, 4123, 4321\end{matrix}$\\
\hline $4_1$ & $1234$ & all~other    \\
\hline $L2a1 \sqcup O$ & $1243$ & $\begin{matrix} 1324, 1342, 2413, 2431,\\ 3124, 3142, 4213, 4231\end{matrix}$   \\
\hline $L2a1 \# L2a1$ & $1243$ & all~other    \\
\hline $L2a1 \sqcup O$ & $2341$ & all  \\
\hline $L6a5$ &$2143$ & $\begin{matrix} 1324, 1342, 2413, 2431,\\ 3124, 3142, 4213, 4231\end{matrix}$ \\
\hline $O \sqcup O\sqcup O \sqcup O\sqcup O$ & $4321$  & all\\
\hline
\end{tabular}
\vskip .2cm
\caption{Links of flat plumbing basket number $4$.} \label{t1}
\end{table}

\begin{thm} \label{class}
\begin{enumerate}
\item[{\rm (1)}]  A link $L$ has the flat plumbing basket number $0$ if and only if $L$ is the trivial knot.

\item[{\rm (2)}]  A link $L$ has the flat plumbing basket number $1$ if and only if $L$ is the trivial link of two components.

\item[{\rm (3)}]  A link $L$ has the flat plumbing basket number $2$ if and only if $L$ is the trivial link of three components.

\item[{\rm (4)}]  A link $L$ has the flat plumbing basket number $3$ if and only if $L$ is either the trivial link of four components or
the Hopf link which is denoted by $L2a1$.

\item[{\rm (5)}] A link $L$ has the flat plumbing basket number $4$ if and only if $K$ is
either the trefoil knot, the figure eight knot, $L2a1 \sqcup O$, $L2a1 \# L2a1$, or $L6a5$.
\end{enumerate}
\begin{proof}
(1) The flat plumbing basket surface with $0$ flat plumbing must be a $2$-disc.
Thus one can easily see that a link $L$ has the flat plumbing basket number $0$ if and only if $L$ is the trivial knot.

(2) The flat plumbing basket surface with $1$ flat plumbing must be a $2$-disc with a single flat annulus.
A link $L$ has the flat plumbing basket number $1$ if and only if $L$ is the trivial link of two components.

(3) Consider a link of the flat plumbing basket number $2$. There are only two possible such surfaces and
one's boundary is the trivial knot and the others is the trivial link of three components. Therefore, we have that
a link $L$ has the flat plumbing basket number $2$ if and only if $L$ is the trivial link of three components.

(4) By considering all possible flat $3$ band diagrams which are $36$ cases,
we find they are either the trivial link of two components, the trivial link of four components or
the Hopf link. However, the trivial link of two components has the flat plumbing basket number $1$.

(5) Suppose a link $L$ has the flat plumbing basket number $4$.
Before we consider all possible permutation presentation $(a_1 a_2 a_3 a_4 : b_1 b_2 b_3 b_4)$,
we first found that all $(a_1 a_2 a_3 a_4 : b_1 b_2 b_3 b_4)$ are reducible except
$(1234:b_1 b_2 b_3 b_4)$, $(1243:b_1 b_2 b_3 b_4)$, $(1342:b_1 b_2 b_3 b_4)$, $(2143: b_1 b_2 b_3 b_4)$,
$(2341:b_1 b_2 b_3 b_4)$ and $(4321:b_1 b_2 b_3 b_4)$. Let us remark that
there are more irreducible permutation presentation but
$(1324:b_1 b_2 b_3 b_4)$ and $(2134:b_1 b_2 b_3 b_4)$ can be obtained from $(1243:b_1 b_2 b_3 b_4)$,
$(1423:b_1 b_2 b_3 b_4)$ and $(3124:b_1 b_2 b_3 b_4)$ can be obtained from $(1423:b_1 b_2 b_3 b_4)$,
$(4123:b_1 b_2 b_3 b_4)$ can be obtained from $(2341:b_1 b_2 b_3 b_4)$ by cyclic relabeling.
Then we find a complete list of links which can be obtained by permutation presentations
of flat plumbing basket surfaces in Table~\ref{t1}. Let us remark that this link corresponding to
permutation presentations $(1342:b_1 b_2 b_3 b_4)$ are all unknot.
Therefore, it complete the proof.
\end{proof}
\end{thm}

Although, we have found the classification theorem of links
of the flat plumbing basket number $4$ or less, it can be used
to determine the flat plumbing basket number of a link which
is either $5$ or $6$.
In fact, if one find
a flat plumbing basket surface of $5$ annuli whose boundary
is a link $L$ which is not listed in the classification theorem,
then the flat plumbing basket number of the link $L$ must be $5$. The following
is a such a example.

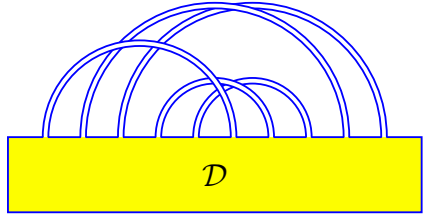
\begin{figure}
$$
\begin{pspicture}[shift=-1.2](-.7,-1.9)(5.2,2)
\psarc[doubleline=true](2.75,0){.75}{-5}{185}
\psarc[doubleline=true](2.25,0){.75}{-5}{185}
\psarc[doubleline=true](2.75,0){1.75}{-5}{185}
\psarc[doubleline=true](2.25,0){1.75}{-5}{185}
\psarc[doubleline=true](1.25,0){1.25}{-5}{185}
\psframe[linecolor=lightgray,fillstyle=solid,fillcolor=lightgray](-.5,-1)(5,0)
\psline(-.03,0)(-.5,0)(-.5,-1)(5,-1)(5,0)(4.53,0)
\psline(.03,0)(.47,0) \psline(.53,0)(.97,0) \psline(1.03,0)(1.47,0)
\psline(1.53,0)(1.97,0) \psline(2.03,0)(2.47,0) \psline(2.53,0)(2.97,0)
\psline(3.03,0)(3.47,0) \psline(3.53,0)(3.97,0) \psline(4.03,0)(4.47,0)
\rput(2.25,-.5){{$\mathcal{D}$}}
\end{pspicture}
$$
\caption{A flat plumbing basket surface of the link $L4a1$ with $5$ flat plumbings.} \label{412fig}
\end{figure}

\begin{cor} \label{link42}
The flat plumbing basket number of the link $L4a1$ is $5$.
\begin{proof}
Since $L4a1$ is not listed in Theorem~\ref{class} and
it is the boundary of the flat plumbing basket surface
$(1,2,3,4,5,1,4,5,2,3)$ as depicted in Figure~\ref{412fig}. It completes the proof
of the corollary.
\end{proof}
\end{cor}

In the case of a knot, the flat plumbing basket number of a knot must be even because
the boundary of a flat plumbing basket surface of $n$ annuli has at most $n+1$ components,
and the number of components is always congruent to $n+1$ modulo $2$.

\begin{cor} \label{52cor}
The flat plumbing basket number of the knot $5_2$ is $6$.
\begin{proof}
The knot $5_2$ is not listed in Theorem~\ref{class} and it admits
a flat plumbing basket surface of $6$ annuli as illustrated in Figure~\ref{52complete}.
Therefore, the flat plumbing basket number of the knot $5_2$ is $6$.
\end{proof}
\end{cor}

Let us remark that the flat plumbing basket number $6$ of the knot $5_2$ in Corollary~\ref{52cor}
is independently found by Hirose and Nakashima~\cite{HN} using the lower bound by the genus
and Alexander polynomial and by Y. Choi, Y. Do and D. Kim~\cite{CDK} using
a complete classification of knots of flat plumbing basket presentation with
$6$ annuli.

\section*{Acknowledgments}

The \TeX\, macro package PSTricks~\cite{PSTricks} was essential for
typesetting the equations and figures. This work was supported by Kyonggi University Research Grant 2011.

\end{document}